\documentclass[preprint,review,12pt]{elsarticle}
\usepackage[margin=2.5cm]{geometry}
\usepackage{hyperref}
\usepackage{amsmath}
\usepackage{lineno}
\usepackage{amsfonts}
\usepackage{amssymb}
\usepackage{bm}
\usepackage{color}
\journal{Journal}
\usepackage{algorithm}
\usepackage{algorithmic}
\usepackage{booktabs}
\usepackage{appendix}

\begin{document}

\begin{frontmatter}

\title{Encoded Forward Backward Stochastic Neural Network for High-Dimensional Backward Stochastic Differential Equations and Parabolic Partial Differential Equations}

%\title{Encoded Forward Backward Stochastic Neural Network for Efficient Computation of Backward Stochastic Differential Equations and High-Dimensional Semilinear Parabolic Partial Differential Equations}

\address[1]{Research Centre for Mathematics and Interdisciplinary Sciences, Shandong University, Qingdao, Shandong Province, 266237, China}
\address[2]{Frontiers Science Center for Nonlinear Expectations, Minister of Education, Shandong University, Qingdao, Shandong Province, 266237, China}

\author[1,2]{Zhao Zhang\corref{cor1}}
\ead{zhaozhang@sdu.edu.cn}
\author[1,2]{Zhuopeng Hou\corref{cor1}}
\ead{202317081@mail.sdu.edu.cn}
\cortext[cor1]{Corresponding author}

\begin{abstract}
Backward stochastic differential equation (BSDE) provides  probabilistic solutions for a class of parabolic partial differential equations (PDEs). DeepBSDE and FBSNN are two deep learning approaches for solving high-dimensional PDEs through approximating the solution of BSDEs. The conventional approach for learning functions defined on continuous domains is via fully-connected networks (FCNs) such that each input dimension is represented by a single neuron. In the current study, a new encoded FBSNN algorithm is proposed to enhance the efficiency and accuracy of approximating BSDEs using encoding and convolution. The input coordinates are encoded as tensors treated as images with multiple channels which can be processed efficiently by convolutional neural networks. The encoding mechanism enriches the input features such that the spatial and temporal features can be balanced. The encoded FBSNN algorithm provides a simple yet effective extension of the vanilla FBSNN algorithm such that BSDEs can be approximated more efficiently. The new algorithm is validated using the essentially high-dimensional Black-Scholes-Barenblatt and Hamilton–Jacobi–Bellman benchmark cases.
\end{abstract}


\end{frontmatter}

%\linenumbers
\section{Introduction}
A backward stochastic differential equation (BSDE) is a type of stochastic differential equation (SDE) with a specified terminal condition, for which the solution must be adapted to an underlying filtration and satisfy both the equation and terminal condition \cite{pardoux1990adapted, peng1993backward, el1997backward}. BSDEs naturally arise in various applications such as stochastic control, mathematical finance and statistical physics \cite{el2008backward}. 
Further, BSDEs provide probabilistic solutions to semilinear and nonlinear parabolic partial differential equations (PDEs) according to the nonlinear Feynman-Kac formula \cite{pardoux1999forward, cheridito2007second}.

In recent years, artificial intelligence and deep learning have achieved remarkable success in many fields including natural language and image processing \cite{lecun2015deep}. Inspired by these achievements, there have been many studies on solving PDEs using deep learning. Among the earliest studies, DeepBSDE algorithms have been proposed to solve high-dimensional parabolic PDEs reformulated as BSDEs, overcoming the curse of dimensionality which is a challenge for conventional numerical methods \cite{han2017deep, han2018solving, beck2019machine}. Almost at the same time, the physics-informed neural network (PINN) method has been proposed to solve PDEs \cite{raissi2017physics, Raissi2017PhysicsID, raissi2019physics}. The basic idea of PINN is to treat neural network (NN) output as PDE solution which is optimized by minimizing the PDE residual evaluated by automatic differentiation (AD). The benefit of PINN is that the learned NN is essentially a meshless surrogate in infinite temporal-spatial domain and can be evaluated at any location. PINN has also been used to solve BSDEs for approximating high-dimensional parabolic PDEs stochastically  without labeled data in the forward backward stochastic neural network (FBSNN) method \cite{raissi2024forward}. 

The convergence of DeepBSDE for decoupled FBSDE systems has been analyzed in \cite{han2020convergence} which is further extended to coupled BSDE systems with jumps in \cite{wang2025deep}. The convergence of FBSNN has not yet been theoretically analyzed according to our literature review. For high-dimensional Hamilton-Jacobi-Bellman equation, a martingale neural network method is proposed based on stochastic optimal control \cite{cai2025soc}. For very-high-dimensional PDEs, AD may be computationally intensive which can be replaced by random finite difference for higher efficiency \cite{cai2025deep}.

In the current study, we focus on decoupled FBSDE systems and propose a new encoded FBSNN algorithm to compute efficiently BSDEs and high-dimensional semilinear parabolic PDEs. Vanilla DeepBSDE and FBSNN algorithms are based on fully connected networks (FCNs). In the new algorithm, the spatial and temporal inputs are encoded into tensors or matrices of multiple channels such that convolutional neural networks (CNNs) can be used to approximate the solutions of BSDEs more efficiently than FCNs. In addition, encoding a scalar temporal coordinate into a vector or matrix enriches the input representation for CNN to learn. Encoding is integrated into the NN model as a layer such that AD can be conducted via back-propagation. Once trained, the CNN-based model can be used to predict the BSDE solution of new sample trajectories of the given SDE.

This paper is organized as follows. First, BSDE and the nonlinear Feynman-Kac formula are introduced as the theoretical basis. Second, the discretization scheme and FBSNN algorithm are reviewed. Third, the encoded FBSNN algorithm is presented, along with a discussion of related works. Finally, the new algorithm is validated on essentially high-dimensional benchmark cases including the Black-Scholes-Barenblatt equation and the Hamilton-Jacobi-Bellman equation.

\section{Methodology}
\subsection{BSDE and the Nonlinear Feynman-Kac Formula}
There have been many studies regarding the link between nonlinear parabolic PDEs and BSDEs in literature \cite{pardoux1990adapted, peng1993backward, el1997backward}. Let $(\Omega, \mathcal{F}, \mathbb{P})$ be a probability space, $\{W_t\}_{t\in[0,T]}$ be a d-dimensional Brownian motion, and $\{\mathcal{F}_t\}_{t\in[0,T]}$ be a filtration generated by $\{W_t\}_{t\in[0,T]}$.
In the current study, we consider decoupled forward-backward SDE (FBSDE) system as
\begin{eqnarray}
    X_t &=&X_0+\int_0^t\mu(s,X_s)ds+\int_0^t\sigma(s,X_s)dW_s \label{4}\\
    Y_t &=&g(X_T)+\int_t^Tf(s,X_s,Y_s,Z_s)ds-\int_t^TZ_sdW_s \label{yt}
\end{eqnarray}
where $t$ is time, $x$ is $d$-dimensional spatial variable, $\mu$ is a known vector-valued function, $\sigma$ is a known $d \times d$ matrix-valued function, $f$ is a known nonlinear function, $\{X_t\}_{t\in[0,T]}$ is a d-dimensional stochastic process, $(Y_t, Z_t)$ is the $\{\mathcal{F}_t\}_{t\in[0,T]}$-adapted solution process. The FBSDE system is related to a class of semilinear PDE written as
\begin{align}
\nonumber
\frac{\partial u}{\partial t}(t,x) + \frac{1}{2}\text{Tr}&\left(\sigma\sigma^T(t,x)(\nabla^2 u)(t,x)\right) + \nabla u(t,x)\cdot \mu(t,x)\\&+ f(t,x,u(t,x),\sigma^T(t,x)\nabla u(t,x)) = 0 \label{2}  
\end{align}
with terminal condition $u(T, x) = g(x)$, where $\text{Tr}$ denotes the matrix trace operator, $\nabla u$ and $\nabla^2 u$ denote the gradient and Hessian of $u$ with respect to $x$, respectively. Under suitable assumptions on the regularity of $\mu$, $\sigma$ and $f$, it holds $\mathbb{P}$-a.s. that
\begin{eqnarray}
    Y_t &=& u(t, X_t) \nonumber\\
    Z_t &=& \nabla u(t,X_t)\sigma(t, X_t)
\end{eqnarray}
which is the nonlinear Feynman-Kac formula \cite{peng1992nonlinear}. Therefore, Eq.~\eqref{yt} can be rewritten as
\begin{align}
\nonumber    
    u(t,X_t) - u(0,X_0)=& - \int_0^t f\left(s,X_s,u(s,X_s),\sigma^T(s,X_s)\nabla u(s,X_s)\right)ds\\&+ \int_0^t [\nabla u(s,X_s)]^T \sigma(s,X_s) dW_s\label{3}.
\end{align}

\subsection{Temporal Discretization}
Given a partition of the time interval $[0, T]$ as $0 = t_0 < t_1 < \cdots < t_N = T$, the Euler-Maruyama scheme is adopted for discretization such that
\begin{align}
    X_{n+1}=X_{n}+\mu(t,X_{n})\Delta t+\sigma(t,X_{n})\Delta W_t ~,\label{8}
\end{align}
where the subscript $n$ denotes time step, $\Delta t=t_{n+1}-t_n$ and $\Delta W_t=W_{t+\Delta t}-W_{t}\sim \mathcal{N}(0,
\Delta t)$. The paths $\{X_{n}\}_{0\leq n\leq N}$ are generated by sampling from Eq.\eqref{8}. For $u(t,X_t)$, Eq.\eqref{3} can be discretized in time accordingly as
\begin{align}
\nonumber
   u(t_{n+1},X_{n+1})=u(t_n,&X_{n})-f\left(t_n,X_{n},u(t_n,X_{n}),\sigma^{T}(t_n,X_{n})\nabla u(t_n,X_{n})\right)\Delta t_n\\&+[\nabla u(t_n,X_{n})]^T\sigma(t_n,X_{n})\Delta W_n \label{13}~.
\end{align}

\subsection{FBSNN Algorithm}
Given the Euler-Maruyama discretization scheme in Eqs.~\eqref{8} and \eqref{13} with terminal condition $u(T, X_T)=g(X_T)$, the idea of FBSNN is to use a NN denoted as $u_\Theta(t, X_t)$ to approximate the mapping $(t, X_{t})\rightarrow u(t,X_{t})$ such that $u_\Theta(t, X_t)$ satisfies Eq.~\eqref{13} and the terminal condition, where $\Theta$ denotes the NN parameters. The NN output $u_\Theta(t,X_{t})$ is substituted into the discrete scheme Eq.\eqref{13} to yield the residual at $(t_n, X_{n})$ as
\begin{align}
\nonumber
   R_{n}=&u_\Theta(t_{n},X_{n})-u_\Theta(t_{n-1},X_{n-1})-{[\nabla u_\Theta(t_{n-1},X_{n-1})]}^T\sigma(t_{n-1},X_{n-1})\Delta W_{n-1}\\&+f(t_{n-1},X_{n-1},u_\Theta(t_{n-1},X_{n-1}),\sigma^T(t_{n-1},X_{n-1})\nabla u_\Theta(t_{n-1},X_{n-1}))\Delta t_{n-1}\label{14}~. 
\end{align}
The BSDE-informed loss functional is given by
\begin{equation}
   \mathcal{L}_{BSDE}= \mathbb{E}\left[\sum_{n=1}^{N}|R_{n}|^2\right]~.\label{15}
\end{equation}
The loss functional for the terminal condition is given by
\begin{align}
  \mathcal{L}_T=\mathbb{E}\left[|u(T,X_T)-g(X_T)|^2\right]\label{19} ~. 
\end{align}
The loss functional for the terminal gradient condition is given by
\begin{align}
    \mathcal{L}_G=\mathbb{E}\left[|\nabla u(T,X_T)-\nabla g(X_T)|^2\right]
\end{align}
The total loss functional is given by
\begin{align}
  \mathcal{L}=\alpha_1\mathcal{L}_{BSDE}+\alpha_2\mathcal{L}_T+\alpha_3\mathcal{L}_G \label{20} ~,
\end{align}
where $\alpha_1$, $\alpha_2$ and $\alpha_3$ are weights for the three loss terms. They are all set to be unity in the current study.

The loss functional is evaluated on a number of sample trajectories simulated using the Euler-Maruyama scheme in Eq.~\eqref{8}. Let $m$ denote a sample trajectory, the BSDE-informed loss functional is evaluated as
\begin{eqnarray}
\nonumber    
    \mathcal{L}\approx \alpha_1\sum_{m=1}^{M}\sum_{n=1}^{N}&|R^m_{n}|^2+\alpha_2\sum_{m=1}^{M}|u(T,X^m_T)-g(X^m_T)|^2\\&+\alpha_3\sum_{m=1}^{M}|\nabla u(T,X^m_T)-\nabla g(X^m_T)|^2
\end{eqnarray}
given $M$ sample trajectories and $N$ time snapshots, where $R^m_n$ is computed at $(t_n, X^m_{n})$ as
\begin{align}
\nonumber
   R_{n}^m&=u_\Theta(t_{n},X^m_{n})-u(t_{n-1},X^m_{n-1})-{[\nabla u(t_{n-1},X^m_{n-1})]}^T\sigma(t_{n-1},X^m_{n-1})\Delta W^m_{n-1}\\&+f(t_{n-1},X^m_{n-1},u(t_{n-1},X^m_{n-1}),\sigma^T(t_{n-1},X^m_{n-1})\nabla u(t_{n-1},X^m_{n-1}))\Delta t_{n-1}~,
\end{align}
where $\nabla u(t_{n-1},X^m_{n-1})$ is calculated explicitly by AD, and $\Delta W^m_{n-1}$ is known from the numerical simulation of Eq.~\eqref{8}. 

\subsection{Encoded FBSNN Algorithm}
DeepBSDE and FBSNN algorithms typically employ FCNs to learn the mapping $(t, X_{t})\rightarrow u_\Theta(t,X_{t})$. In this study, we encode $(t, X_{t})$ into tensor forms which can be treated as images with multiple channels such that CNN-based neural structures can be employed to enhance the convergence rate of training. 
Compared to FCN, CNN-based structures reduce the number of trainable parameters by  local connectivity and weight sharing. This not only mitigates overfitting but also simplifies the optimization landscape for enhanced convergence rate. In addition, CNN-based structures inherently build a hierarchical representation to capture multiscale features for data with spatial coherence. 

The spatial input $X_t$ and temporal coordinate $t$ are encoded separately. For $t$, the sine-cosine positional encoding technique is applied to enrich the time-series features. For each snapshot $t$, a temporal encoding vector $TE(t)$ is generated by
\begin{eqnarray}
    TE_{2i}(t) &=& t+sin(\frac{t}{10000^{\frac{2i}{d}}})\nonumber\\
    TE_{2i+1}(t) &=& t+cos(\frac{t}{10000^{\frac{2i}{d}}})\label{code}
\end{eqnarray}
where $d$ denotes the encoding dimension, $i=1,2,...,\frac{d}{2}$ is the encoding-dimension index. The vector $TE(t)$ is then reshaped into a matrix.

\begin{figure}[htbp]    
    \centering     
    \includegraphics[width=0.99\textwidth]{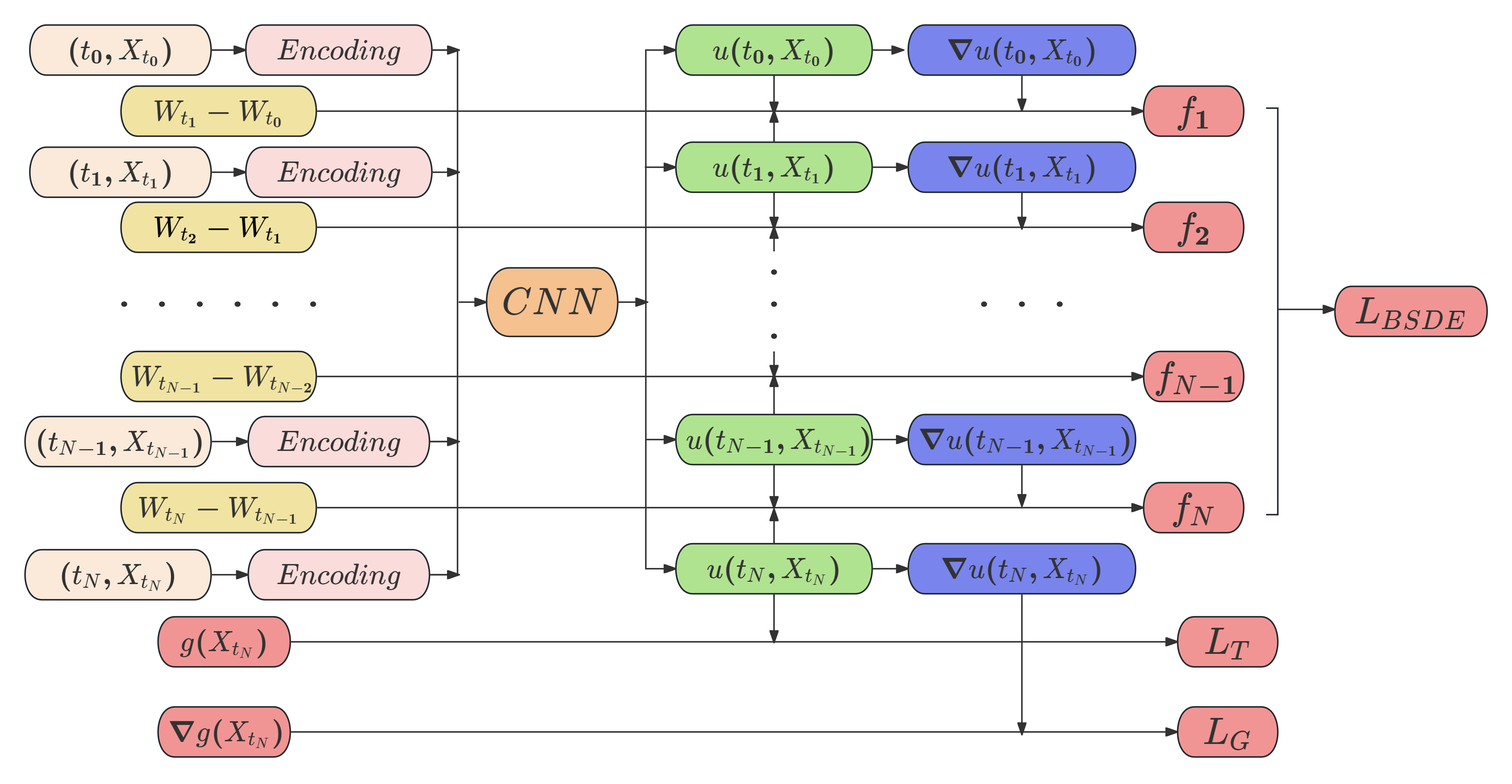} 
    \caption{Illustration of the encoded FBSNN structure. The time domain is discretized into $N$ intervals, with each row corresponding to a specific time step $t=t_0,t_1,...,t_{N-1},t_N$. }     
    \label{one} 
\end{figure}

For the spatial input $X_t$ that is high-dimensional, encoding each dimension as in Eq.~\eqref{code} leads to an image with a large number of channels which is less efficient. In fact, as $X_t$ is a vector itself, we can simply reshape $X_t$ into an image-lime matrix. In order to increase the dimension of the encoded matrix for $X_t$, linear interpolation can be used to obtain a matrix of higher dimensions containing richer features for CNNs to learn. The matrices of $TE(t)$ and encoded $X_t$ are concatenated as a tensor representing an image of two channels to be the input of CNN.
The Adam optimizer is adopted for minimizing the loss functional \cite{adam2014method}. The encoded FBSNN algorithm is presented in Algorithm \ref{alg-cnn} and illustrated in Fig.~\ref{one}.

\begin{algorithm}
	\caption{Encoded FBSNN Algorithm for decoupled FBSDE}
	\label{alg-cnn}
	\begin{algorithmic}[1]
		\REQUIRE Number of time steps $N$ for discretizing $[0, T]$, number of sample trajectories $M$, number of epochs $N_{epochs}$, the initial point $X_0$ and the terminal condition $g(x)$.
		\ENSURE CNN model $\mathcal{N}_{\Theta}$ for approximating the nonlinear mapping from encoded $(t, X_{t})$ to $u(t,X_{t})$.
        \STATE Generate $M$ sample trajectories according to the Euler-Maruyama discretization scheme in Eq.\eqref{8} for the training set.
        \STATE Substitute the terminal state $X_N$ into the given terminal condition $g(x)$ and its gradient $\nabla g(x)$ to obtain the corresponding labels. 
        \FOR{epoch $j=1$ to $N_{epochs}$}
        \STATE Given $t_0$ and $X_0$, generate the encoded tensor and then predict $Y_0$ by $\mathcal{N}_{\Theta}$. Use AD to compute $Z_0$.
        \STATE Retrieve a mini-batch of samples from the training set.
		\FOR{time step $i = 1$ to $N$}
        \STATE Encode $t_n$ and $X_n$ into a tensor.
        \STATE Predict $Y_i$ by $\mathcal{N}_{\Theta}$ and obtain $Z_i$ by AD.
		\STATE Compute the BSDE loss $R_i$ for time step $t_i$. 
		\ENDFOR
		\STATE Accumulate the BSDE losses $R_i$ of all time steps to obtain the total BSDE loss.
        \STATE Calculate the terminal loss $\mathcal{L}_T$  and gradient loss $\mathcal{L}_G$ for the predicted values of $Y_N$ and $Z_N$ using the corresponding labels. 
        \STATE Calculate the total loss function as $\mathcal{L}=\alpha_1\mathcal{L}_{BSDE}+\alpha_2\mathcal{L}_T+\alpha_3\mathcal{L}_G$.
        \STATE Use Backpropagation to calculate the gradient of total loss $\mathcal{L}$ with respect to the neural network parameters $\Theta$. 
        \STATE Use Adam optimizer to update $\Theta$ and minimize the total loss $\mathcal{L}$.
        \ENDFOR
	\end{algorithmic}
\end{algorithm}

In encoded FBSNN, the mapping $(t, X_{t})\rightarrow u_\Theta(t,X_{t})$ is approximated by
\begin{equation}
    u_\Theta(t,X_{t})=\mathcal{N}_{\Theta}\circ\mathcal{E}_{encoder}(t,X_{t})
\end{equation}
where $\mathcal{N}_{\Theta}$ denote the CNN,  $\Theta$ denotes trainable parameters, and $\mathcal{E}_{encoder}$ is the encoding block. Although $\mathcal{E}_{encoder}$ contains no trainable parameters, it is included in the computational graph such that back-propagation can be conducted.
For $\mathcal{N}_{\Theta}$, a compact structure containing an output block and two convolutional blocks is built as
\begin{align}
    \mathcal{N}_{\Theta}:= \mathcal{L}_{out}\circ \mathcal{B}_{2}\circ \mathcal{B}_{1}
\end{align}
where $\mathcal{L}_{out}$ is the output block and $\mathcal{B}_i$, $i=1,2$, denote the convolutional blocks such that
\begin{eqnarray}
    \mathcal{B}_1 &:=& \mathcal{P}_1\circ \Phi \circ \mathcal{BN}\circ Conv2d\nonumber\\
    \mathcal{B}_2 &:=& \mathcal{P}_2\circ \Phi \circ \mathcal{BN}\circ Conv2d\nonumber \\
    \mathcal{L}_{out} &:=& Linear\circ \Phi \circ Linear\nonumber\\
\end{eqnarray}
where $\Phi$ denotes the activation function, $\mathcal{BN}$ is batch normalization, $Linear$ means linear layers, $\mathcal{P}_1$ and $\mathcal{P}_2$ are upsampling operators such as max-pooling and average-pooling. The details of the CNN structure for encoded FBSNN are presented in Fig.~\ref{net} in \ref{a1}. The details of the FCN structure in the vanilla FBSNN algorithm for comparison are presented in Fig.~\ref{fcn} in \ref{a1}.

\subsection{Comparison of Related Work}
In DeepBSDE \cite{han2018solving}, the first-order-gradient term $[\nabla u(t_n,X_{n})]^T\sigma(t_n,X_{n})$ is not obtained via AD but predicted by the NN, while the loss functional only accounts for the misfit of terminal condition such that the value of a single snapshot for $u(t, X_t)$ is obtained after training. If the value of $u$ at a different snapshot is desired, the training of NN needs to be conducted again. Further, the number of sub-networks grows linearly with the number of time steps.

In FBSNN \cite{raissi2024forward}, the gradient of $u(t_n,X_{n})$ is obtained by AD. Both the governing BSDE equation and terminal condition are accounted for in the loss functional. After training, each $u(t_n,X_{n})$ for $0\le n \le N$ is obtained. All time steps share the same network such that the network size can stay unchanged for increased number of time steps. The vanilla FBSNN is based on FCN which represents each input dimension using a single neuron.

In encoded FBSNN, the temporal and spatial inputs are encoded into matrices such that the input features can be enriched. The matrices for temporal and spatial inputs are of the same dimension to balance the learning of temporal and spatial features. The matrices are concatenated into a tensor treated as an image of two channels which can be learned efficiently by CNNs benefiting from local connectivity and weight sharing.

\section{Benchmark Cases}
The new encoded FBSNN algorithm is validated on two essentially high-dimensional PDE benchmark cases reported in \cite{han2018solving}.
\subsection{Black-Scholes-Barenblatt Equation}
\label{test1}
\subsubsection{Model}
The Black-Scholes-Barenblatt (BSB) equation\cite{volatilitiespricing} is a robust extension of the classical Black-Scholes model. It provides a theoretical foundation for the emergence of bid-ask spreads. In the presence of model uncertainty, such as unknown volatility that may vary within a certain range, the BSB equation offers a worst-case protection strategy for option pricing and hedging. A core principle of this model is that traders price options under the assumption that market parameters would be adversely selected. To hedge this risk, convex positions are priced using the upper bound of volatility, while concave positions are priced using the lower bound. This approach represents a paradigm shift in financial engineering, moving from the pursuit of precise pricing towards robust hedging \cite{meyer2006black}.

The BSB equation in $[0, T]\times R^{100}$ can be written as
\begin{align}
   \frac{\partial u}{\partial t}&(t,x)+r·x·\nabla u(t,x)+\frac{{\sigma}^2}{2}\sum_{i=1}^d|x_i|^2\frac{\partial^2u}{\partial x_i^2}(t,x)-r·u(t,x)=0\label{21} ~.
\end{align} 
where $T=1,~\sigma=0.4, ~r=0.05,~\xi=(1,0.5,1,0.5,...,1,0.5)\in R^{100},~g(x)={||x||}^2$, and the terminal condition $u(T,x)=g(x)$.
The corresponding FBSDE is written as
\begin{eqnarray}
    dX_t &=& \sigma X_tdW_t,~t\in [0,T] \nonumber\\
    X_0 &=& \xi\nonumber\\
    dY_t &=& r(Y_t-\frac{Z_t}{\sigma}X_t)dt+Z_tX_tdW_t,~t\in[0,T]\nonumber\\
    Y_T &=& g(X_T)
\end{eqnarray}
This equation admits the explicit solution
\begin{equation}
u(t,x)=exp((r+{\sigma}^2)(T-t))g(x)
\end{equation}
which is used as the ground truth to validate the accuracy of machine learning algorithm.

\subsubsection{Settings}
\label{test1set}
For validation, the new encoded FBSNN algorithm is compared to the vanilla FBSNN method presented in \cite{raissi2024forward}. 
For training, there are in total 50 time steps and 5000 sample trajectories. In each epoch, 100 sample trajectories are randomly selected to form a minibatch. A two-stage learning rate is employed. A higher learning rate is used in the initial phase to achieve rapid convergence, followed by a lower learning rate for fine-tuning. 

The relative error between the approximated $u_\Theta(t,x)$ and ground truth on the test set is adopted to quantify the prediction accuracy. A sufficiently low test error guarantees that the mapping $(t, X_{t})\rightarrow u(t,X_{t})$ can be approximated accurately by $u_\Theta(t,x)$.
It is worth noting that the training and testing datasets are kept unchanged for different algorithms to ensure fair comparison. 

The encoded input dimension is set to be $2\times 20\times 20$ such that there are 2 channels for spatial and temporal coordinates. The 100-dimensional spatial coordinate is reshaped into $10\times 10$ and then bilinearly interpolated into $20\times 20$. The scalar temporal coordinate is encoded by the sine-cosine approach in Eq.~\eqref{code} to be a 400-dimensional vector which is reshaped into $20\times 20$.

The learning rate is set as follows for optimal prediction accuracy. For the encoded FBSNN, the learning rate of the first 2000 epochs is 0.001 which is reduced to be 0.0001 subsequently. For the vanilla FBSNN, the learning rate of the first 10000 epochs is 0.001 which is reduced to be 0.0001 subsequently. The default settings of hyperparameters are as in Table~\ref{tab1}.

\begin{table}
    \centering
    \scalebox{0.8}{
    \begin{tabular}{|c|c|c|}
    \hline
         & Encoded FBSNN & Vanilla FBSNN\\
         \hline
        Train set size & 5000 & 5000\\
         \hline
         Minibatch size & 100 & 100\\
         \hline
        Test set size & 1000 & 1000\\
         \hline
        Epochs & 3000 & 12000\\
         \hline
        1st Learning rate & 0.001 for 2000 epochs & 0.001 for 10000 epochs \\
         \hline
         2nd Learning rate & 0.0001 for remaining epochs & 0.0001 for remaining epochs \\
         \hline
    \end{tabular}
    }
    \caption{Settings for Encoded and vanilla FBSNN in the Black-Scholes-Barenblatt test case.}
    \label{tab1}
\end{table}

\subsubsection{Results and Discussion}
\begin{figure}[htbp]    
    \centering     
    \includegraphics[width=0.9\textwidth]{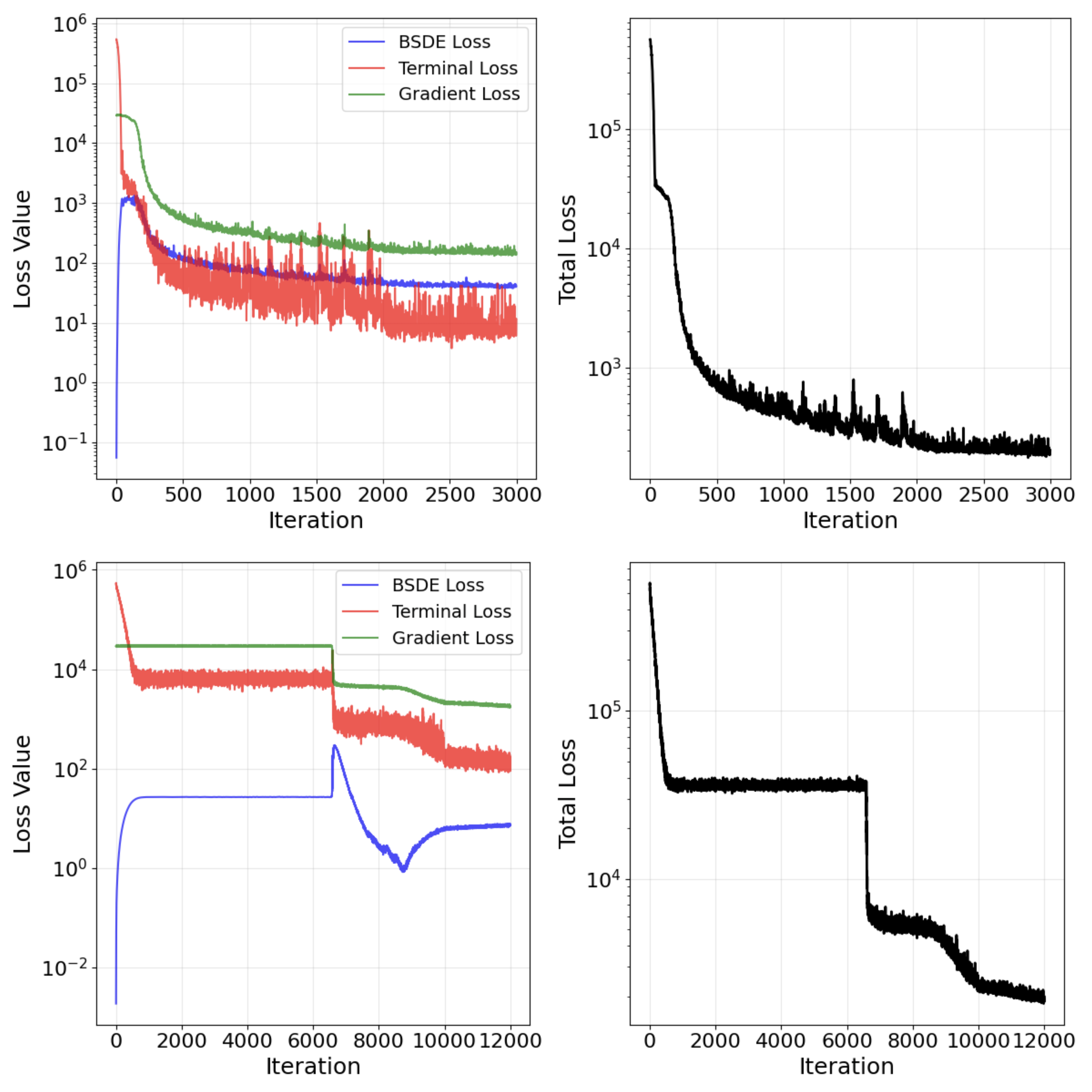} 
    \caption{The convergence of $\mathcal{L}_{BSDE},\mathcal{L}_{T}$ and $\mathcal{L}_{G}$ over training epochs for the encoded FBSNN algorithm (top left) and the vanilla FBSNN algorithm (bottom left) for the Black-Scholes-Barenblatt test case. The convergence of total loss $\mathcal{L}$ with respect to training epochs for the encoded FBSNN algorithm (top right) and the vanilla FBSNN algorithm (bottom right). The encoding dimension is $20\times 20$.}     
    \label{two} 
\end{figure}

The convergence of the training loss with respect to the number of epochs for encoded and vanilla FBSNN algorithms are illustrated in Fig.~\ref{two}. The encoded algorithm requires less than about one third of epochs to reach a lower training loss than the vanilla algorithm. For the training of both algorithms, the terminal loss is decreasing monotonically, while the BSDE loss increases first and then drops to a relatively lower level. This implies that the initial solution in a local minima satisfies the BSDE approximately but not the terminal condition. Nonetheless, the total loss drops monotonically driven by the stochastic Adam optimizer.

The results of encoded and vanilla FBSNN with different training strategies are summarized in Table~\ref{tab:bs_errors}. To account for the fluctuations of stochastic optimization using mini batches, the results are averaged over 10 independent runs. In the table, encoded FBSNN is the algorithm proposed in the current study. The CNN structure for encoded FBSNN is detailed in Fig.~\ref{net}. Reshaped FBSNN means that the 100-dimensional spatial input is simply reshaped into a $10\times 10$ matrix rather than being encoded into a $20\times 20$ matrix while the other structures are kept as the same as encoded FBSNN. The training and prediction accuracies of reshaped FBSNN are both lower than that of encoded FBSNN. The encoded FBSNN  algorithm can achieve more than 40\% increase in prediction accuracy with less computational cost compared to vanilla FBSNN. 

\begin{table}[htbp]
  \centering
  \scalebox{0.8}{
  \begin{tabular}{lcccc}
    \toprule
    Method & Epochs & Relative error(train) & Relative error(test)  & Cost \\
    \midrule
    Vanilla FBSNN & \(12000\) & \( 0.98\%\) & \( 1.03\% \) & 21.4 min \\
    Reshaped FBSNN & \(3000\)  & \( 1.11\%\) & \( 1.16\% \) & 12.6 min \\
    Encoded FBSNN &\(3000\) & \(0.58\% \) & \( 0.61\%\)  & 16.2 min \\
    \bottomrule
  \end{tabular}
  }
    \caption{Summary of results using encoded and vanilla FBSNN algorithms with different settings in the Black-Scholes-Barenblatt test case. Encoded FBSNN is the algorithm proposed in the current study. Reshaped FBSNN means that the spatial input is simply reshaped into a  matrix rather than encoded.}
    \label{tab:bs_errors}
\end{table}

\begin{figure}[htbp]    
    \centering     
    \includegraphics[width=0.99\textwidth]{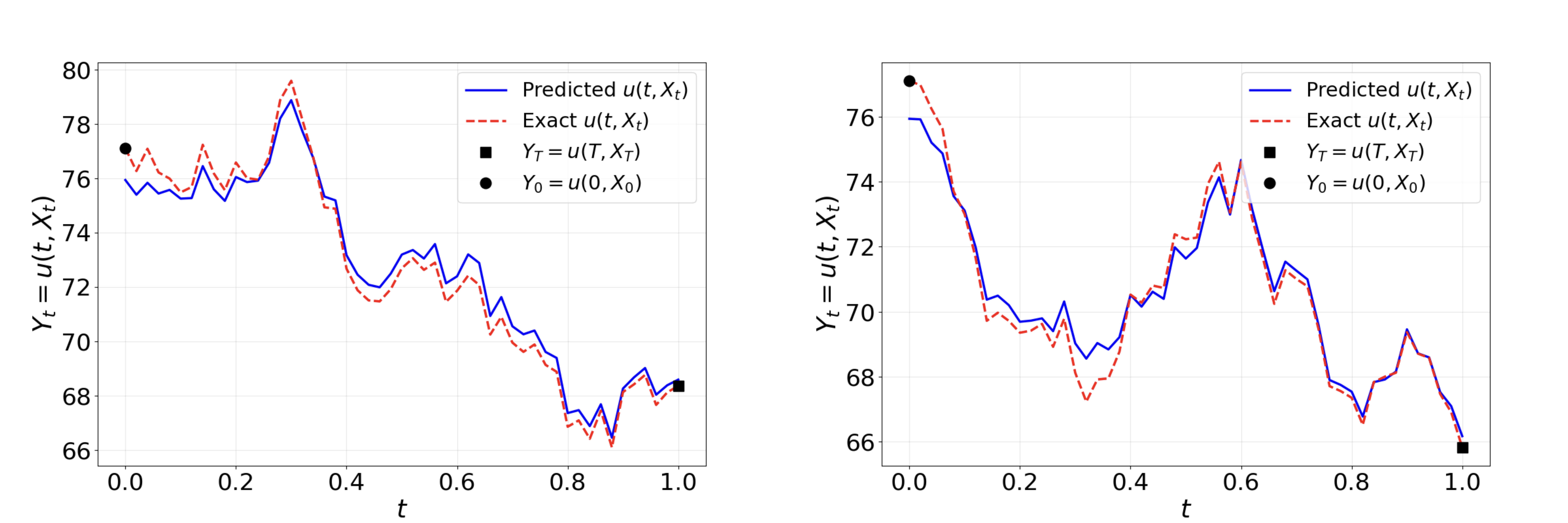} 
    \caption{Representative predicted sample trajectory for $Y_t=u(t,X_t)$ using encoded FBSNN in the Black-Scholes-Barenblatt test case. One learned trajectory of medium accuracy from the training set (left), and one predicted trajectory of medium accuracy from the test set (right) are plotted against the ground truth.}     
    \label{three} 
\end{figure}

To examine visually the training and prediction accuracy of BSDE trajectories approximated by encoded FBSNN, one trajectory of medium accuracy is selected from each of the training and test sets to be plotted against the explicit solution in Fig.~\ref{three}. The approximated trajectories are in good agreement with the exact solution. 

\subsection{Hamilton–Jacobi–Bellman Equation}
\subsubsection{Model}
In the stochastic optimal control theory, the Hamilton-Jacobi-Bellman (HJB) equation provides a rigorous framework for deriving globally optimal closed-loop feedback policies in dynamic environments with uncertainty. It enables the synthesis of real-time, state-dependent control laws that guarantee optimal performance.
Considering a classical linear-quadratic Gaussian (LQG) control problem in 100 dimensional space,
\begin{align}
  dX_t=2\sqrt{\lambda}m_tdt+\sqrt{2}dW_t\label{23}  ,
\end{align}
with $t\in[0,T],X_0=x$ and the cost function 
\begin{align}
   J(\{m_t\}_{0 \leq t \leq T}) = \mathbb{E} \left[ \int_0^T \|m_t\|^2 \, dt + g(X_T) \right]\label{24} ~,
\end{align}
where
$\{X_t\}_{t \in [0,T]}$ is the state process,
$\{m_t\}_{t \in [0,T]}$ is the control process,
$\lambda$ is a positive constant representing the strength of the control and
$\{W_t\}_{t \in [0,T]}$ is a standard Brownian motion.
The target is to minimize the cost functional through the control process. This test case has been used in \cite{han2018solving}.

The HJB equation for this problem is given by
\begin{align}
 \frac{\partial u}{\partial t}(t,x) + \Delta u(t,x) - \lambda||\nabla u(t,x)||^2 = 0 \label{25} , 
\end{align}
with the terminal condition $u(T,x)=g(x)=\ln\left(\frac{1+\|x\|^2}{2}\right)$. The value of the solution $u(t, x )$ at $t = 0$ represents the optimal cost when the state starts from $x$ \cite{yong1999stochastic}. According to the Feyman-Kac formula, one can show that the analytical solution of Eq.\eqref{25} admits the explicit formula
\begin{align}
  u(t,x) = -\frac{1}{\lambda}\ln\left(E\left[\exp\left(-\lambda g\left(x+\sqrt{2}W_{T-t}\right)\right)\right]\right)\label{26} ~.
\end{align}
Monte Carlo simulation using the analytical solution is used to provide the reference solution for comparison.

\subsubsection{Settings}
For validation, the new encoded FBSNN algorithm is compared to the vanilla FBSNN algorithm. 
Similarly as in Section \ref{test1}, there are in total 50 time steps. A two-stage learning rate is employed. A higher learning rate is used in the initial phase to achieve rapid convergence, followed by a lower learning rate for fine-tuning. The input temporal and spatial coordinates are encoded into a $2\times 20\times 20$ tensor similarly as in Section \ref{test1set} to be learned by a CNN.
The number of epochs, learning rate, training and test set size are set as in Table~\ref{tab1}.

\subsubsection{Results and Discussion}
\begin{figure}[htbp]    
    \centering     
    \includegraphics[width=0.9\textwidth]{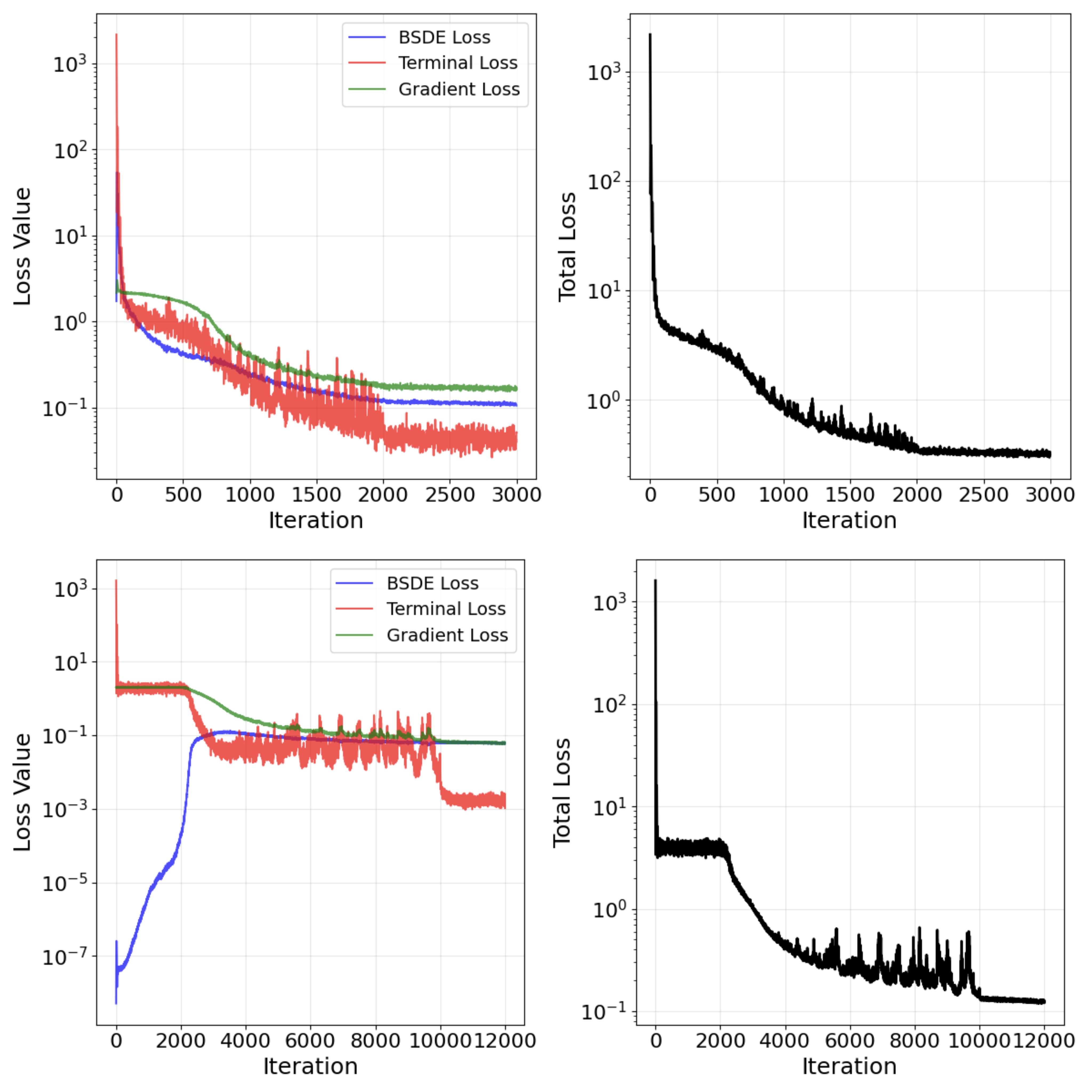} 
    \caption{The convergence of $\mathcal{L}_{BSDE},\mathcal{L}_{T}$ and $\mathcal{L}_{G}$ over training epochs for the encoded FBSNN algorithm (top left) and the vanilla FBSNN algorithm (bottom left) for the Hamilton–Jacobi–Bellman test case. The convergence of total loss $\mathcal{L}$ with respect to training epochs for the encoded FBSNN algorithm (top right) and the vanilla FBSNN algorithm (bottom right). The encoding dimension is $20\times 20$.}     
    \label{five} 
\end{figure}

The convergence of the training loss over epochs for encoded and vanilla FBSNN are shown in Fig.~\ref{five}. With fewer training iterations and less computational time, the encoded FBSNN algorithm achieves a much lower loss than the vanilla FBSNN algorithm. Although the total and terminal losses generally present monotone decrease despite small-scale fluctuations, the BSDE losses generally increases first before dropping to lower levels. 

The results of encoded and vanilla FBSNN with different settings are summarized in Table \ref{tab:hjb_errors}. To account for the fluctuations of stochastic optimization using mini batches, the results are averaged over 10 independent runs.
The encoded FBSNN algorithm achieves higher accuracy compared to the vanilla FBSNN algorithm using much fewer training epochs and less computational cost. The training and prediction accuracies of reshaped FBSNN are both lower than that of encoded FBSNN. The encoded FBSNN algorithm can achieve much lower relative errors representing more than 100\% increase in accuracy using less epochs and computational cost compared to the vanilla FBSNN algorithm.

\begin{table}[htbp]
  \centering
  \scalebox{0.8}{
  \begin{tabular}{lcccc}
    \toprule
    Method & Epochs & Relative error (train) & Relative error (test)  & Cost \\
    \midrule
    Vanilla FBSNN & \(12000\) & \( 0.42\% \) & \( 0.44\% \)  & 21.4 min \\
    Reshaped FBSNN & \(3000\) & \( 0.34\% \) & \( 0.36\% \)  & 12.6 min \\
    Encoded FBSNN & \(3000\) & \( 0.20\% \) & \( 0.21\% \)  & 16.2 min \\
    \bottomrule
  \end{tabular}
  }
  \caption{Summary of results using encoded and vanilla FBSNN algorithms with different settings in the Hamilton–Jacobi–Bellman test case. Encoded FBSNN is the algorithm proposed in the current study. Reshaped FBSNN means that the spatial input is simply reshaped into a  matrix rather than encoded.}
  \label{tab:hjb_errors}
\end{table}

\begin{figure}[htbp]    
    \centering     
    \includegraphics[width=0.99\textwidth]{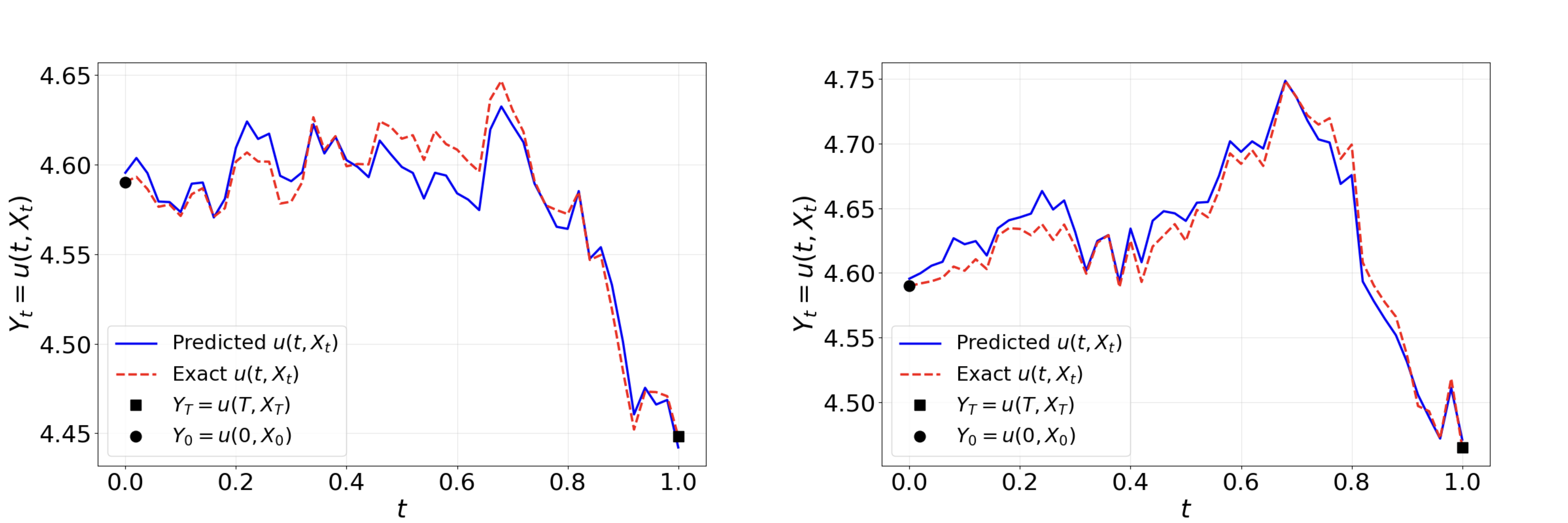} 
    \caption{Representative predicted sample trajectory for $Y_t=u(t,X_t)$ using encoded FBSNN in the Hamilton–Jacobi–Bellman test case. One learned trajectory of medium accuracy from the training set (left), and one predicted trajectory of medium accuracy from the test set (right) are plotted against the ground truth.}     
    \label{six} 
\end{figure}

To examine BSDE trajectories approximated by encoded FBSNN, one trajectory of medium accuracy is selected from each of the training and test sets to be plotted against the explicit solution in Fig.~\ref{six}. The approximated trajectories are in good agreement with the exact solution. 

\section{Conclusions}
In this paper, a new encoded FBSNN algorithm has been proposed to solve high-dimensional BSDEs and the corresponding semilinear parabolic PDEs over the entire spatiotemporal domain. Compared to the vanilla FBSNN method, the new algorithm encodes the input spatial-temporal coordinates into tensors that can be treated as images such that CNN-based neural networks can be used to approximate BSDE solutions efficiently on continuous domains. The features of the input coordinates can be enriched by encoding and CNN is generally more efficient than FCN for image-like inputs. The effectiveness of encoded FBSNN for enhancing efficiency and accuracy in solving high-dimensional BSDEs and semilinear parabolic PDEs have been validated on benchmark cases of Black-Scholes-Barenblatt and Hamilton-Jacobi-Bellman equations. The encoded FBSNN algorithm demonstrates a simple yet effective extension of the vanilla FBSNN method such that modern CNN-based network structures can be employed for solving BSDEs efficiently.

\section*{Code Availability}
The encoded FBSNN codes are available at "https://github.com/felix0426/Encoded-FBSNN". 

\section*{Acknowledgements}
The research is supported by the Natural Science Foundation of Shandong Province (No.ZR2024MA057), the Fundamental Research Funds for the Central Universities and the Future Plan for Young Scholars of Shandong University.

\bibliography{ref}

@article{pardoux1990adapted,
  title={Adapted solution of a backward stochastic differential equation},
  author={Pardoux, Etienne and Peng, Shige},
  journal={Systems \& control letters},
  volume={14},
  number={1},
  pages={55--61},
  year={1990},
  publisher={Elsevier}
}

@article{raissi2019physics,
	title={Physics-informed neural networks: A deep learning framework for solving forward and inverse problems involving nonlinear partial differential equations},
	author={Raissi, Maziar and Perdikaris, Paris and Karniadakis, George E},
	journal={Journal of Computational Physics},
	volume={378},
	pages={686--707},
	year={2019},
	publisher={Elsevier}
}

@article{peng1993backward,
  title={Backward stochastic differential equations and applications to optimal control},
  author={Peng, Shige},
  journal={Applied Mathematics and Optimization},
  volume={27},
  number={2},
  pages={125--144},
  year={1993},
  publisher={Springer}
}

@article{el1997backward,
  title={Backward stochastic differential equations in finance},
  author={El Karoui, Nicole and Peng, Shige and Quenez, Marie Claire},
  journal={Mathematical finance},
  volume={7},
  number={1},
  pages={1--71},
  year={1997},
  publisher={Wiley Online Library}
}

@misc{el2008backward,
  title={Backward stochastic differential equations and applications},
  author={El Karoui, Nicole and Hamad{\`e}ne, Said and Matoussi, Anis},
  journal={Indifference pricing: theory and applications},
  volume={8},
  pages={267--320},
  year={2008},
  publisher={Springer Berlin}
}

@article{cheridito2007second,
  title={Second-order backward stochastic differential equations and fully nonlinear parabolic PDEs},
  author={Cheridito, Patrick and Soner, H Mete and Touzi, Nizar and Victoir, Nicolas},
  journal={Communications on Pure and Applied Mathematics},
  volume={60},
  number={7},
  pages={1081--1110},
  year={2007},
  publisher={Wiley Online Library}
}

@article{pardoux1999forward,
  title={Forward-backward stochastic differential equations and quasilinear parabolic PDEs},
  author={Pardoux, Etienne and Tang, Shanjian},
  journal={Probability theory and related fields},
  volume={114},
  number={2},
  pages={123--150},
  year={1999},
  publisher={Springer}
}

@inproceedings{peng1992nonlinear,
  title={A nonlinear Feynman-Kac formula and applications},
  author={Peng, Shige and others},
  booktitle={Proceedings of Symposium of System Sciences and Control Theory},
  pages={173--184},
  year={1992},
  organization={World Scientific}
}

@article{lecun2015deep,
  title={Deep learning},
  author={LeCun, Yann and Bengio, Yoshua and Hinton, Geoffrey},
  journal={nature},
  volume={521},
  number={7553},
  pages={436--444},
  year={2015},
  publisher={Nature Publishing Group UK London}
}

@article{beck2019machine,
  title={Machine learning approximation algorithms for high-dimensional fully nonlinear partial differential equations and second-order backward stochastic differential equations},
  author={Beck, Christian and E, Weinan and Jentzen, Arnulf},
  journal={Journal of Nonlinear Science},
  volume={29},
  number={4},
  pages={1563--1619},
  year={2019},
  publisher={Springer}
}

@article{han2017deep,
  title={Deep learning-based numerical methods for high-dimensional parabolic partial differential equations and backward stochastic differential equations},
  author={Han, Jiequn and Jentzen, Arnulf and others},
  journal={Communications in mathematics and statistics},
  volume={5},
  number={4},
  pages={349--380},
  year={2017},
  publisher={Springer}
}

@article{han2018solving,
  title={Solving high-dimensional partial differential equations using deep learning},
  author={Han, Jiequn and Jentzen, Arnulf and E, Weinan},
  journal={Proceedings of the National Academy of Sciences},
  volume={115},
  number={34},
  pages={8505--8510},
  year={2018},
  publisher={National Academy of Sciences}
}

@article{raissi2017physics,
  title={Physics informed deep learning (part i): Data-driven solutions of nonlinear partial differential equations},
  author={Raissi, Maziar and Perdikaris, Paris and Karniadakis, George Em},
  journal={arXiv preprint arXiv:1711.10561},
  year={2017}
}

@article{Raissi2017PhysicsID,
  title={Physics Informed Deep Learning (Part II): Data-driven Discovery of Nonlinear Partial Differential Equations},
  author={Maziar Raissi and Paris Perdikaris and George Em Karniadakis},
  journal={ArXiv},
  year={2017},
}

@article{cai2025deep,
  title={Deep random difference method for high dimensional quasilinear parabolic partial differential equations},
  author={Cai, Wei and Fang, Shuixin and Zhou, Tao},
  journal={arXiv preprint arXiv:2506.20308},
  year={2025}
}

@article{cai2025soc,
  title={SOC-MartNet: A Martingale Neural Network for the Hamilton--Jacobi--Bellman Equation Without Explicit in Stochastic Optimal Controls},
  author={Cai, Wei and Fang, Shuixin and Zhou, Tao},
  journal={SIAM Journal on Scientific Computing},
  volume={47},
  number={4},
  pages={C795--C819},
  year={2025},
  publisher={SIAM}
}

@article{wang2025deep,
  title={Deep Learning Numerical Methods for High-Dimensional Quasilinear PIDEs and Coupled FBSDEs with Jumps},
  author={Wang, Wansheng and Wang, Jie and Li, Jinping and Gao, Feifei and Fu, Yi and Ye, Zaijun},
  journal={SIAM Journal on Scientific Computing},
  volume={47},
  number={3},
  pages={C706--C737},
  year={2025},
  publisher={SIAM}
}

@article{han2020convergence,
  title={Convergence of the deep BSDE method for coupled FBSDEs},
  author={Han, Jiequn and Long, Jihao},
  journal={Probability, Uncertainty and Quantitative Risk},
  volume={5},
  number={1},
  pages={5},
  year={2020},
  publisher={Springer}
}

@incollection{raissi2024forward,
  title={Forward--backward stochastic neural networks: deep learning of high-dimensional partial differential equations},
  author={Raissi, Maziar},
  booktitle={Peter Carr Gedenkschrift: Research Advances in Mathematical Finance},
  pages={637--655},
  year={2024},
  publisher={World Scientific}
}

@article{adam2014method,
  title={A method for stochastic optimization},
  author={Adam, Kingma DP Ba J and others},
  journal={arXiv preprint arXiv:1412.6980},
  volume={1412},
  number={6},
  year={2014}
}

@book{yong1999stochastic,
  title={Stochastic controls: Hamiltonian systems and HJB equations},
  author={Yong, Jiongmin and Zhou, Xun Yu},
  volume={43},
  year={1999},
  publisher={Springer Science \& Business Media}
}

@article{volatilitiespricing,
  title={PRICING AND HEDGING DERIVATIVE SECURITIES IN MARKETS WITH UNCERTAIN VOLATILITIES},
  author={VOLATILITIES, UNCERTAIN}
}

@article{meyer2006black,
  title={The Black Scholes Barenblatt equation for options with uncertain volatility and its application to static hedging},
  author={Meyer, Gunter H},
  journal={International Journal of Theoretical and Applied Finance},
  volume={9},
  number={05},
  pages={673--703},
  year={2006},
  publisher={World Scientific}
}

%\newpage
\appendix

\section{Neural Network Structures}
\label{a1}

\begin{figure}[htbp]    
    \centering     
    \includegraphics[width=0.73\textwidth]{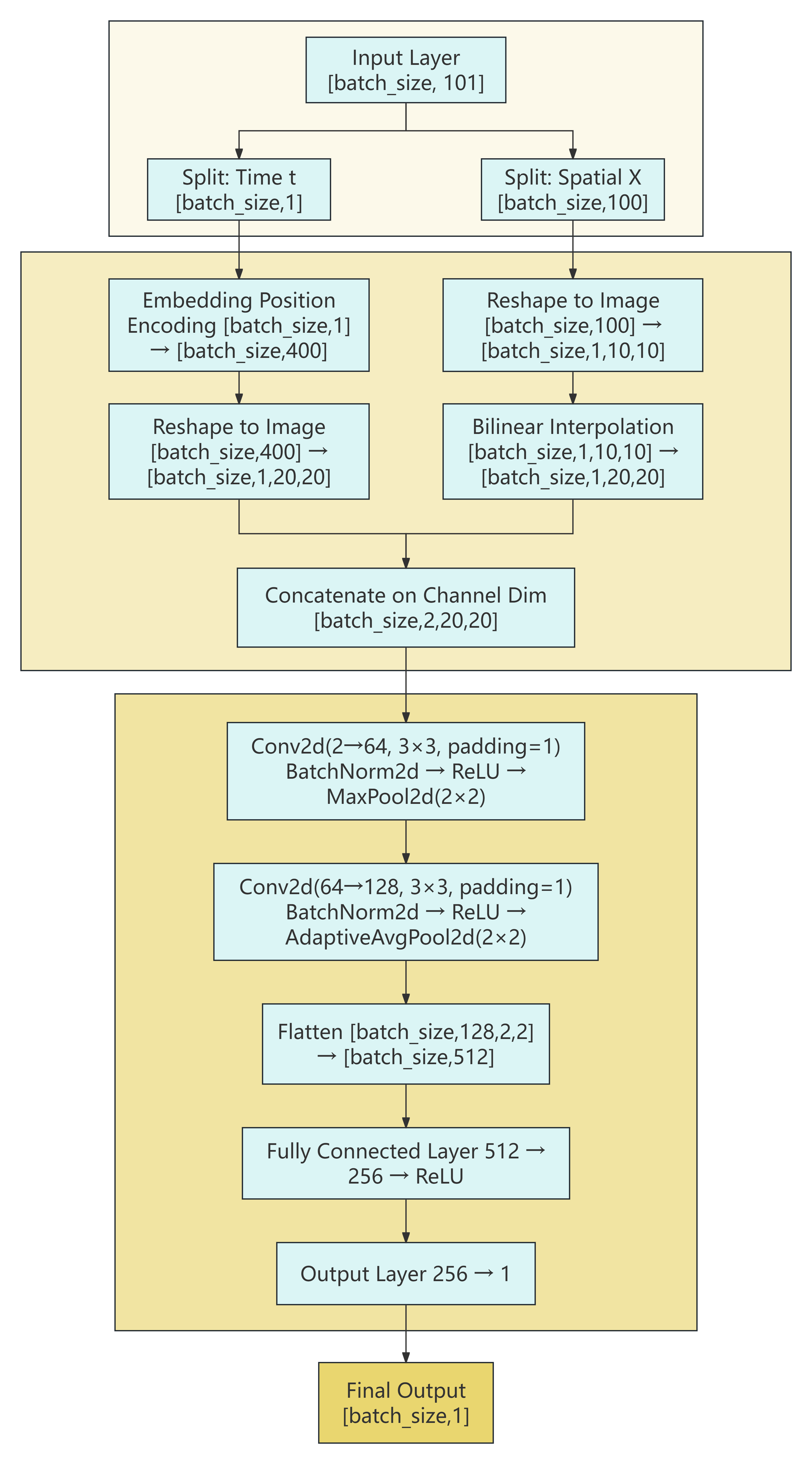} 
    \caption{Illustration of the encoding block and the CNN network structure. The second upsampling using AdaptiveAvgPool2d can be replaced by MaxPool2d if necessary.}  
    \label{net} 
\end{figure}

Fig.\ref{net} illustrates the network architecture for for the case where the encoding dimension is set to $20\times 20$.
First, the initial convolutional layer expands the input image of two channels to a feature map of 64 channels using $3\times3$ kernels, followed by batch normalization and nonlinear activation using ReLU function. The map is then downsampled via $2\times2$ max-pooling to a lower resolution. Next, the second convolutional layer further increases the channel count to 128, followed by batch normalization and ReLU activation. The feature map is subsequently compressed into a uniform $2\times2$ grid using adaptive average pooling. Finally, the flattened 512-dimensional feature vector is fed into a fully connected layer with 256 neurons and ReLU activation, followed by a linear layer to be output as a scalar. 

\begin{figure}[htbp]    
    \centering     
    \includegraphics[width=0.43\textwidth]{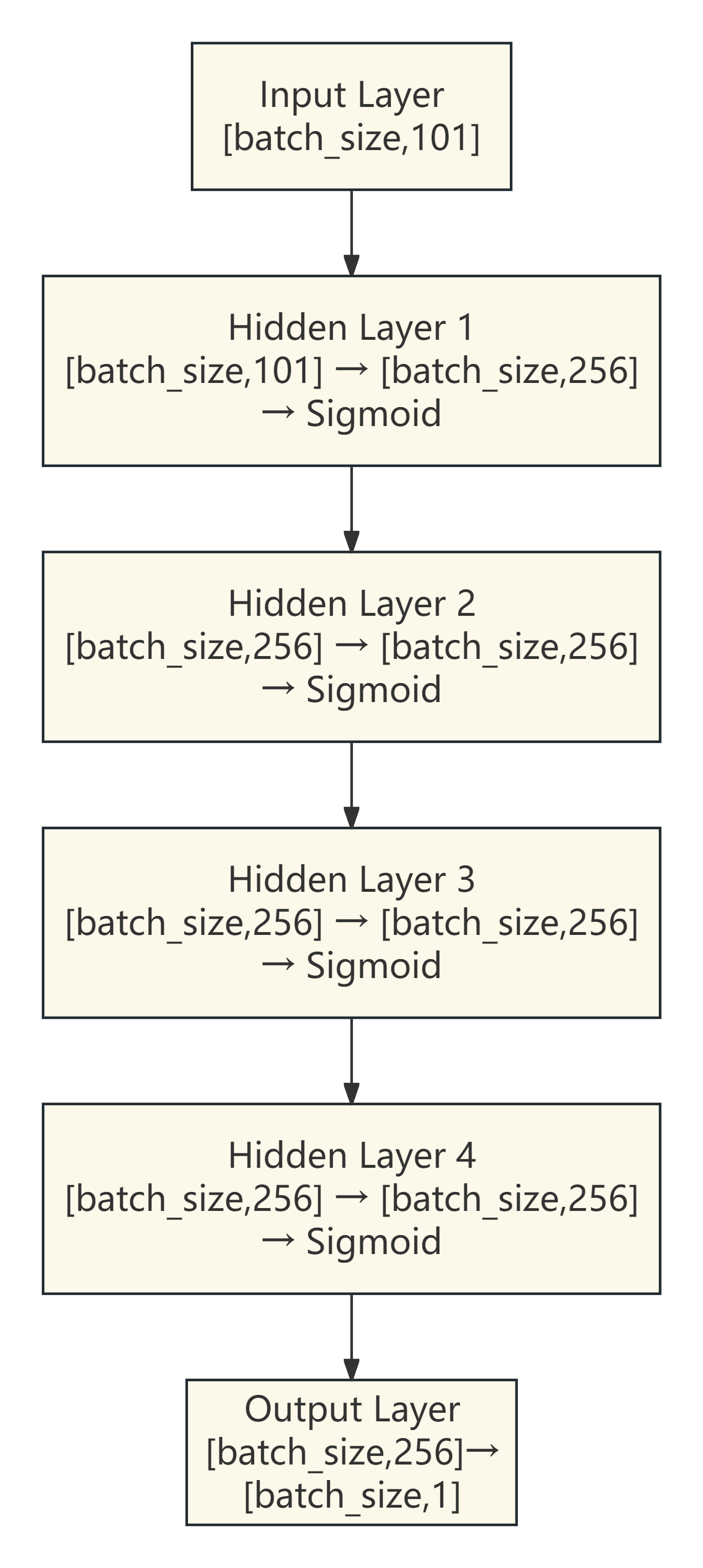} 
    \caption{Illustration of the FCN structure employed in the vanilla FBSNN algorithm for comparative study. Other activation functions can be used if necessary.}  
    \label{fcn} 
\end{figure}

Fig.\ref{fcn} illustrates the FCN structure employed in the vanilla FBSNN Algorithm. The FCN consists of 6 layers. The input is the concatenation of the one-dimensional temporal variable t and the 100-dimensional spatial variable X. Each of the four hidden layers consists of 256 neurons. The comparative results in the current study are obtained using Sigmoid activation function. Tanh and ReLU activation functions can be used but the results are approximately unchanged. 

\end{document}